\documentclass[12pt,sort&compress]{article}
\usepackage{amsfonts}
\usepackage{bm}
\usepackage{makecell}
\usepackage{amssymb,amsmath}
\usepackage{amscd}
\usepackage{lipsum} 
\usepackage{latexsym}
\usepackage{CJK}
\usepackage{longtable}
\def\g{\gamma}

\def\cl{\centerline}

\def\vs{\vspace*}

\def\Z{\mathbb{Z}}
\def\CC{\mathbb{C}}
\def\C{\mathbb{C}}

\def\QED{\hfill$\Box$}

\def\t{\tilde}

\def\pa{\partial}
\def\la{\lambda}
\def\SV{\mathcal{SV}}

\numberwithin{equation}{section}
\newtheorem{theo}{Theorem}[section]
\newtheorem{defi}[theo]{Definition}
\newtheorem{coro}[theo]{Corollary}
\newtheorem{lemm}[theo]{Lemma}

\newtheorem{rema}[theo]{Remark}

\newtheorem{example}[theo]{Example}
\newtheorem{proposition}[theo]{Proposition}
\makeatletter
\newcommand\blfootnote[1]{% 
\begingroup 
\renewcommand\thefootnote{}\footnote{#1}% 
\addtocounter{footnote}{-1}% 
\endgroup 
}

\def\@biblabel#1{#1.~}

\makeatother

\begin{document}
\vs{10pt} \cl{\large {\bf Cohomology of %a class of rank two
the Schr\"odinger-Virasoro conformal algebra}}
%\footnote{Corresponding author: Henan Wu (wuhenan@sxu.edu.cn).}} \vs{12pt}Cohomology of the Schr\"odinger-Virasoro conformal algebra
\cl{% Lamei Yuan$^{\,\ddag}$,
Henan Wu$^{1}$, Lipeng Luo$^{2}$}
% \cl{\small{ $^{\ddag}$ Academy of Fundamental and Interdisciplinary
% Sciences,}}\cl{\small{Harbin Institute of Technology, Harbin 150080, China}}
\cl{\small{$^{1}$School of Mathematical Sciences, Shanxi University, Taiyuan, 030006, China}}
\cl{\small{$^{2}$School of Mathematical Sciences, Tongji University, Shanghai, 200092, China}}
\cl{\small E-mail: wuhenan@sxu.edu.cn, luolipeng@tongji.edu.cn
 }\vs{6pt}
\small\parskip .005 truein \baselineskip 3pt \lineskip 3pt

\noindent{\bf Abstract:} Both the basic cohomology groups and the reduced cohomology groups of the Schr\"odinger-Virasoro conformal algebra with trivial coefficients are completely determined.
\blfootnote{The second author is the corresponding author: Lipeng Luo (luolipeng@tongji.edu.cn)}

\vs{5pt}

\vs{5pt}

\noindent{\bf MR(2000) Subject Classification:}~ 17B65, 17B68

%\vs{18pt}
\section{Introduction}
%%%%%%%%%%%%%%%
Lie conformal algebra, which was introduced by Kac in \cite{K1,K2}, gives an axiomatic description of the operator product expansion (or rather its Fourier transform) of chiral fields in conformal field theory (see \cite{BPZ}). It has been shown that the theory of Lie conformal algebras has close connections to vertex algebras, infinite-dimensional Lie algebras satisfying the locality property in \cite{KacL} such as affine Kac-Moody algebras, the Virasoro algebra, and Hamiltonian formalism in the theory of nonlinear evolution equations (see \cite{BDK,DK,K2,SYX, SY,Z1,Z2}).
The general structure theory, representation theory of Lie conformal algebras were systematically developed in \cite{CK,DK}.

In this paper, we focus on the cohomology theory of Lie conformal algebras. The general cohomology theory of conformal algebras with coefficients in an arbitrary conformal module was
developed in \cite{BKV}, where explicit computations of cohomology groups for the Virasoro conformal algebra and current conformal algebra were given. Some low dimensional cohomology groups of the general Lie
conformal algebras $gc_N$ were studied in \cite{S}.
All cohomology groups of the Heisenberg-Virasoro conformal algebra with trivial coefficients were determined in \cite{YW}. It was also shown in \cite{BKV} that the basic cohomology groups of Lie conformal algebra are naturally isomorphic to those of its annihilation Lie algebra.
However, the arbitrary dimensional cohomology groups of a Lie algebra, especially of a Lie algebra with infinite dimension, are very difficult to compute. So the the cohomology theory of Lie conformal algebra actually gives an efficient method to determine the cohomology groups of certain Lie algebras. In this paper, we will compute the cohomology of the Schr\"odinger-Virasoro conformal algebra(see Definition \ref{def3.1}), which was first introduced in \cite{SY} and has a close  relationship with the well-known Schr\"odinger-Virasoro Lie algebra.
Our methods may be useful to determine higher cohomology groups of some infinite dimensional Lie algebras. This is the main motivation to present our work.

The rest of the paper is organized as follows. In Section 2, we recall some basic definitions, notations, and related known results about Lie conformal algebras. %Then the main theorem of this paper is given. 
In Section 3, we determine the basic cohomology groups of the Schr\"odinger-Virasoro conformal algebra with coefficients in its trivial module $\mathbb{C}_a$. In Section 4, we compute the reduced cohomology groups of the Schr\"odinger-Virasoro conformal algebra with coefficients in its module $\mathbb{C}_a$. As a byproduct, we also compute the reduced cohomology groups with coefficients in $M_{\alpha,\beta}$ in case $\beta\neq0$.
Our main results are summarized in Theorems \ref{th1}, \ref{th2}, \ref{th3}, \ref{th4}.

Throughout this paper, we use notations $\mathbb{C}$, $\mathbb{Z}$ and $\mathbb{Z_{+}}$ to represent the set of complex numbers, integers and nonnegative integers, respectively. In addition, all vector spaces and tensor products are over $\mathbb{C}$. In the absence of ambiguity, we abbreviate $\otimes_{\mathbb{C}}$ to $\otimes$.
%%%%%%%%%%%%%%%
\vs{8pt}

\section{Preliminaries}

   In this section, we recall some basic definitions, notations and related results about Lie conformal algebras for later use. For a detailed description, one can refer to \cite{BKV,YW}. %Then we present our main results.
\subsection{Lie conformal algebra}
\begin{defi}
\begin{em}(\cite{DK}) 
A Lie conformal algebra $\mathcal {A}$ is a $\CC[\partial ]$-module endowed with a $\CC$-bilinear map $\mathcal {A}\otimes \mathcal {A}\rightarrow \CC[\lambda]\otimes \mathcal {A}$, $a\otimes b \mapsto [a_\lambda b]
$ subject to the following relations ($a, b, c\in \mathcal {A}$):
\begin{align*}
[\partial a_\lambda b]&=-\lambda[a_\lambda b],\ \ [ a_\lambda \partial b]=(\partial+\lambda)[a_\lambda b],% \ \ \mbox{(conformal\  sesquilinearity)},\label{Lc1}
\\
{[a_\lambda b]} &= -[b_{-\lambda-\partial}a],% \ \ \mbox{(skew-symmetry)},\label{Lc2}
\\
{[a_\lambda[b_\mu c]]}&=[[a_\lambda b]_{\lambda+\mu
}c]+[b_\mu[a_\lambda c]].%\ \ \mbox{(Jacobi \ identity)}\label{Lc3}.
\end{align*}
\end{em}
\end{defi}

The Lie conformal algebra of our study is the so-called Schr\"odinger-Virasoro conformal algebra introduced in \cite{SY}. And its definition is given in the following.
\begin{defi}\label{def3.1}
\begin{em}
The Schr\"odinger-Virasoro conformal algebra is a finite Lie conformal algebra
$\mathcal{SV}=\mathbb{C}[\partial]L\bigoplus
\mathbb{C}[\partial]M\bigoplus\mathbb{C}[\partial]Y$, endowed with the
following  nontrivial $\lambda$-brackets\vspace*{-7pt}
\begin{eqnarray}
&&[L_\lambda L]=(\partial+2\lambda)L,\ \ \ [Y_\la
L]=(\frac12\pa+\frac32\la)Y, \label{lamda-bracket11}\\
&&{[L_\lambda Y]}=(\partial+\frac32\lambda)Y,\ \ \,{[Y_\lambda
Y]}=(\partial+2\lambda)M,\label{lamda-bracket22}\\[3pt]
&&{[L_\lambda M]}=(\partial+\lambda)M,\ \ \,[M_\la L]=\la
M.\label{lamda-bracket33}
\end{eqnarray}
\end{em}
\end{defi}
The Schr\"odinger-Virasoro conformal algebra $\mathcal{SV}$ contains the
Virasoro conformal algebra $Vir=\mathbb{C}[\partial]L$, whose representation theory and cohomology theory were investigated in \cite{CK} and \cite{BKV}, respectively.
Also, $\mathcal{SV}$ contains the
Heisenberg-Virasoro conformal algebra $\mathcal{HV}=\mathbb{C}[\partial]L\bigoplus
\mathbb{C}[\partial]M$, whose representation theory and cohomology theory were investigated in \cite{WY} and \cite{YW}, respectively.
\subsection{Conformal module}
\begin{defi}
\begin{em}(\cite{CK})
A conformal module $V$ over a Lie conformal algebra $\mathcal {A}$
is a $\mathbb{C}[\partial]$-module equipped with a $\CC$-bilinear map
$\mathcal {A}\otimes V\rightarrow
V[\lambda]$, $a\otimes v\mapsto a_\lambda v$, satisfying the following relations for any $a,b\in\mathcal {A}$, $v\in V$,
\begin{eqnarray*}
&&a_\lambda(b_\mu v)-b_\mu(a_\lambda v)=[a_\lambda b]_{\lambda+\mu}v,\\
&&(\partial a)_\lambda v=-\lambda a_\lambda v,\ a_\lambda(\partial
v)=(\partial+\lambda)a_\lambda v.
\end{eqnarray*}
%If $V$ is finitely generated over $\mathbb{C}[\partial]$, then $V$ is simply called \emph {finite}.
\end{em}
\end{defi}

\begin{example} Let $\mathcal{A}$ be an arbitrary Lie conformal algebra and $a\in\C$.  Then $\mathcal{A}$ admits a family of $1$-dimensional modules $\CC_a$ defined by
$$\CC_a=\CC,\ \partial v=a v,\ \mathcal {A}_\lambda v=0,\ \forall v\in\CC_a.$$
And we abbreviate $\C_0$ to $\C$ in the sequel. It is easy to check that the modules $\mathbb{C}_a$ with $a \in \mathbb{C}$ exhaust all trivial irreducible $\mathcal{A}$-modules.

\end{example}

The classification of all finite irreducible nontrivial $\SV$-modules was obtained in \cite{WY} .

 \begin{proposition}\label{prr}
 All free nontrivial $\SV$-modules of rank one over $\mathbb{C}[\partial]$ are as follows($ \alpha,\beta \in \mathbb{C}$):
  \begin{align*}
M_{\alpha,\beta}=\mathbb{C}[\partial]v,\qquad L_\lambda v=(\partial+\alpha\lambda+\beta)v,\qquad Y_\lambda v=M_\lambda v=0.
 \end{align*}
 Moreover, the module $M_{\alpha,\beta}$ is irreducible if and only if $\alpha$ is non-zero. The module $M_{0,\beta}$ contains a unique nontrivial submodule $ (\partial+\beta)M_{0,\beta}$ isomorphic to $M_{1,\beta}$. The modules $M_{\alpha,\beta}$ with $\alpha\neq 0$ exhaust all finite irreducible nontrivial $\SV$-modules.
\end{proposition}

\subsection{Basic cohomology}
\begin{defi}\label{cochain}\rm (\cite{BKV}) An \emph{$n$-cochain} ($n\in\Z_+$) of a Lie conformal algebra $\mathcal{A}$ with coefficients in an
$\mathcal{A}$-module $V$ is a $\CC$-linear map\vs{-5pt}
\begin{eqnarray*}
\gamma:\mathcal{A}^{\otimes n}\rightarrow V[\la_1,\cdots,\la_n],\ \
\ a_1\otimes\cdots \otimes a_n \mapsto
\g_{\la_1,\cdots,\la_n}(a_1,\cdots,a_n)
\end{eqnarray*}
satisfying the following conditions:\begin{itemize}\parskip-3pt
\item[\rm(1)] $\g_{\la_1,\cdots,\la_n}(a_1,\cdots,\pa a_i,\cdots,
a_n)=-\la_i\g_{\la_1,\cdots,\la_n}(a_1,\cdots, a_n)$ \ (conformal antilinearity),
\item[\rm (2)] $\g$ is skew-symmetric with respect to simultaneous permutations
of $a_i$'s and $\la_i$'s \ (skew-symmetry).
\end{itemize}
\end{defi}

%As usual, let $\mathcal{A}^{\otimes 0}= \CC$, so that a $0$-cochain
%is an element of $V$.
Denote by  ${\t C}^n(\mathcal {A},V)$ the set
of all $n$-cochains. The differential $d_n$ of an $n$-cochain $\g$ is
defined as follows:
\begin{align}\label{ddd}
&(d_n\g)_{\la_1,\cdots,\la_{n+1}}(a_1,\cdots,a_{n+1})\nonumber\\
&=\mbox{$\sum\limits_{i=1}^{n+1}$}(-1)^{i+1}a_{i_{\la_i}}\g_{\la_1,\cdots,\hat{\la_i},\cdots,\la_{n+1}}(a_1,\cdots,\hat{a_i},\cdots,a_{n+1})\nonumber\\
&+\mbox{$\sum\limits_{1\le i<j\le n+1}$}(-1)^{i+j}\g_{\la_i+\la_j,\la_1,\cdots,\hat{\la_i},\cdots,\hat{\la_j},\cdots,\la_{n+1}}([a_{i_{\la_i}}a_j],a_1,\cdots,\hat{a_i},\cdots,\hat{a_j},\cdots,a_{n+1}),
\end{align}
where $\g$ is linearly extended over the polynomials in $\la_i$. %In
%particular, if $\g\in V$ is a $0$-cochain, then
%$(d_0\g)_\la(a)=a_\la\g$.

It was shown in \cite{BKV} that the operator $d$ preserves the space of cochains and $d^2=0$. Thus the cochains of a Lie conformal algebra $\mathcal{A}$ with coefficients in an $\mathcal{A}$-module $V$ form a complex, called the {\it basic complex} :
\begin{equation}
\begin{CD}
\cdot\cdot\cdot@>>> \t C^{n-1}(\mathcal{A},V)  @>{\rm d}_{n-1}>> \t C^{n}(\mathcal{A},V) @>{\rm d}_n>>\t C^{n+1}(\mathcal{A},V)  @>>>\cdot\cdot\cdot
\end{CD}
\end{equation}
%which will be denoted by ${\t C}^\bullet(\mathcal {A},V)=\bigoplus_{n\in \mathbb{Z}_+}{\t C}^n(\mathcal {A},V)$.
The related cohomology is called the {\it basic cohomology} of the Lie conformal algebra $\mathcal{A}$ with coefficients
in its module $V$ and denoted by ${\rm \t H}^q (\mathcal{A},V),\ q\in\Z_+$, more details see the following definition.

\begin{defi}
\begin{em}
An element $\gamma$ in ${\tilde C}^q(\mathcal{A},V)$ is called
a {\it $q$-cocycle} if $d(\gamma)=0$; a {\it $q$-coboundary} if there exists a $(q-1)$-cochain $\phi\in\tilde
C^{q-1}(\mathcal{A},V)$ such that $\gamma=d(\phi)$. Two cochains
$\gamma_1$ and $\gamma_2$ are called {\it equivalent} if $\gamma_1-\gamma_2$ is a coboundary. \end{em}
\end{defi}

Denote by $\tilde D^q(\mathcal{A},V)$ and $\tilde B^q(\mathcal{A},V)$ the spaces
of $q$-cocycles and $q$-boundaries, respectively. Then, we can obtain that
\begin{eqnarray*}
{\rm \tilde H}^q(\mathcal{A},V)=\tilde D^q(\mathcal{A},V)/\tilde B^q(\mathcal{A},V)=\{\mbox{equivalent classes of
$q$-cocycles}\}.
\end{eqnarray*}

\subsection{Reduced cohomology}
Moreover, one can define a (left) $\CC[\partial]$-module structure on $\t
C^n(\mathcal{A},V)$ by
\begin{eqnarray*}
(\partial\g)_{\la_1,\cdots,\la_n}(a_1,\cdots,
a_n)=(\partial_V+\mbox{$\sum\limits_{i=1}^n$}\la_i)\g_{\la_1,\cdots,\la_n}(a_1,\cdots,
a_n),
\end{eqnarray*}
where $\partial_V$ denotes the action of $\partial$ on $V$. Then $d\partial=\partial d$
and %thus $\partial{\t C}^\bullet(\mathcal {A},V)\subset {\t C}^\bullet(\mathcal {A},V)$, i.e.,
\begin{equation}
\begin{CD}
\cdot\cdot\cdot@>>>\partial  \t C^{n-1}(\mathcal{A},V)  @>{\rm d}_{n-1}>>\partial  \t C^{n}(\mathcal{A},V) @>{\rm d}_n>>\partial \t C^{n+1}(\mathcal{A},V)  @>>>\cdot\cdot\cdot
\end{CD}
\end{equation}
forms a subcomplex of the basic complex. The quotient complex %$C^\bullet(\mathcal {A},V)={\t C}^\bullet(\mathcal {A},V)/\partial{\t C}^\bullet(\mathcal {A},V)=\bigoplus_{n\in \mathbb{Z}_+}C^n(\mathcal {A},V)$, i.e.,
\begin{equation}
\begin{CD}
\cdot\cdot\cdot@>>>
\frac{\t C^{n-1}(\mathcal{A},V)}{\partial \t C^{n-1}(\mathcal{A},V)}
@>{\bar{\rm d}}_{n-1}>>
\frac{\t C^{n}(\mathcal{A},V)}{\partial \t C^{n}(\mathcal{A},V)}
@>{\bar{\rm d}}_n>>
\frac{\t C^{n+1}(\mathcal{A},V)}{\partial \t C^{n+1}(\mathcal{A},V)}  @>>>\cdot\cdot\cdot
\end{CD}
\end{equation}
%\begin{eqnarray*}
%\rightarrow\frac{\t C^{n-1}(\mathcal{A},V)}{\partial \t C^{n-1}(\mathcal{A},V)}\rightarrow \frac{\t C^{n}(\mathcal{A},V)}{\partial \t C^{n}(\mathcal{A},V)}\rightarrow\frac{\t C^{n+1}(\mathcal{A},V)}{\partial \t C^{n+1}(\mathcal{A},V)}\rightarrow
%\end{eqnarray*}
is called the {\it reduced complex}. And its cohomology is called the \emph{reduced
cohomology} of the Lie conformal algebra $\mathcal{A}$ with coefficients
in $V$ and denoted by ${\rm H}^q (\mathcal{A},V),\ q\in\Z_+$.

\begin{rema}
The basic cohomology ${\rm \t H}^q (\mathcal{A},V)$ is naturally a $\CC[\partial]$-module, whereas the reduced cohomology ${\rm H}^q (\mathcal{A},V)$ is a complex vector space.
\end{rema}
The exact sequence $0\longrightarrow \partial\t C^\bullet \stackrel{i}{\longrightarrow} \t C^\bullet \stackrel{p}{\longrightarrow}C^\bullet \longrightarrow 0$ gives a long exact sequence of the cohomology groups:
\begin{align}\label{longexact}
0\longrightarrow& H^0(\partial\t C^\bullet) \stackrel{i_0}{\longrightarrow} \t H^0(\mathcal{A},V) \stackrel{p_0}{\longrightarrow} H^0 (\mathcal{A},V) \longrightarrow\nonumber\\
\longrightarrow& H^1(\partial\t C^\bullet) \stackrel{i_1}{\longrightarrow} \t H^1(\mathcal{A},V) \stackrel{p_1}{\longrightarrow} H^1 (\mathcal{A},V) \longrightarrow\\
\longrightarrow& H^2(\partial\t C^\bullet) \stackrel{i_2}{\longrightarrow} \t H^2(\mathcal{A},V) \stackrel{p_2}{\longrightarrow} H^2 (\mathcal{A},V) \longrightarrow\cdots.\nonumber
\end{align}

%Inspired by \cite{CK}, we can obtain the following result, which plays an important role in the classification of the reduced cohomology of  Lie conformal algebra $\mathcal{A}$ with coefficients in $V$.
\begin{proposition}(\cite{BKV})\label{pro2.12}
In degrees $\geq1$, the complexes $\t C^\bullet$ and $\partial\t C^\bullet$ are isomorphic under the map
\begin{align}
\t C^\bullet \to \partial\t C^\bullet, \quad \g \mapsto \partial\cdot \g.
\end{align}
Therefore, $H^q(\partial\t C^\bullet)\cong\t H^q(\mathcal{A},V)$ for $q\geq1$.
\end{proposition}

\begin{rema}
%The above proposition does not imply that
In the long exact sequence (\ref{longexact}),
the maps $H^q(\partial\t C^\bullet)\to\t H^q(\mathcal{A},V)$ induced by the embedding $\partial\t C^\bullet\subset\t C^\bullet$ are not isomorphisms.
\end{rema}

\section{Basic cohomology of $\mathcal{SV}$ with trivial coefficients}
In this section, we will compute the basic cohomology groups of $\mathcal{SV}$ with coefficients in its trivial module $\mathbb{C}_a$.
Since ${\rm \tilde
H}^q(\mathcal{SV},\mathbb{C}_a)\cong {\rm \tilde
H}^q(\mathcal{SV},\mathbb{C})$ for any $a\in \C$ (\cite{BKV}), we only need to compute ${\rm \tilde
H}^q(\mathcal{SV},\mathbb{C})$.
In this case, by (\ref{ddd}), the differential $d_n$ of a $n$-cochain $\g$ is
given as follows:
\begin{eqnarray*}
&&(d_n\g)_{\la_1,\cdots,\la_{n+1}}(a_1,\cdots,a_{n+1})\nonumber\\&&\ \ \ =\mbox{$\sum\limits_{1\le i<j\le n+1}$}(-1)^{i+j}\g_{\la_i+\la_j,\la_1,
\cdots,\hat{\la_i},\cdots,\hat{\la_j},\cdots,\la_{n+1}}([a_{i_{\la_i}}a_j],
a_1,\cdots,\hat{a_i},\cdots,\hat{a_j},\cdots,a_{n+1}).
\end{eqnarray*}

\begin{lemm}\label{l0} ${\rm \tilde
H}^0(\mathcal{SV},\mathbb{C})={\rm H}^0(\mathcal{SV},\mathbb{C})=\mathbb{C}$.
\end{lemm}
{\it Proof.}
For any $\gamma\in \tilde
C^0(\mathcal{SV},\mathbb{C})=\mathbb{C}$, $(d_0\gamma)_\lambda (a)=a_\lambda \gamma =0$ for
$a\in \mathcal{SV}$. This means $\tilde D^0(\mathcal{SV},\mathbb{C})=\mathbb{C}$ and $\tilde B^0(\mathcal{SV},\mathbb{C})=0$. Thus ${\rm \tilde
H}^0(\mathcal{SV},\mathbb{C})=\mathbb{C}$. Moreover, $ {\rm H}^0(\mathcal{SV},\mathbb{C})=\mathbb{C}$
since $\partial\mathbb{C}=0$.
\QED

Let $\gamma\in \tilde
C^q(\mathcal{HV},\mathbb{C})$ with $q>0$. By Definition \ref{cochain}, $\gamma$ is determined by its value on $a_1\otimes\cdots \otimes a_q$ with $a_i\in\{L, Y, M\}$. Since $\gamma$ is skew-symmetric, we can always assume that the first $k$ variables are $L$, the middle $l$ variables are $Y$ and the last $m$ variables are $M$ in $\gamma_{\lambda_1,\cdots,\lambda_{q}}(a_1,\cdots,a_q)$.
Thus we can regard $\gamma_{\lambda_1,\cdots,\lambda_{q}}(a_1,\cdots,a_q)$ as a polynomial in $\lambda_1,\cdots,\lambda_q$, which is skew-symmetric in $\lambda_1,\cdots,\lambda_k$,  in $\lambda_{k+1},\cdots,\lambda_{k+l}$, and in $\lambda_{k+l+1},\cdots,\lambda_{k+l+m}$, respectively, where $q=k+l+m$. Therefore, $\gamma_{\lambda_1,\cdots,\lambda_{q}}(a_1,\cdots,a_q)$ is divisible by $$\mbox{$\prod\limits_{1\leq i< j\leq k}$}(\lambda_i-\lambda_j)\mbox{$\prod\limits_{1\leq i< j\leq l}$}(\lambda_{k+i}-\lambda_{k+j})\mbox{$\prod\limits_{1\leq i< j\leq m}$}(\lambda_{k+l+i}-\lambda_{k+l+j}),$$ whose degree is $\begin{pmatrix}k\\2\end{pmatrix}+\begin{pmatrix}l\\2\end{pmatrix}+\begin{pmatrix}m\\2\end{pmatrix}$.

Following \cite{BKV}, we define an operator $\tau:\t C^q(\mathcal{HV},\CC)\rightarrow
\t C^{q-1}(\mathcal{HV},\CC)$ by
\begin{eqnarray}
(\tau
\g)_{\la_1,\cdots,\la_{q-1}}(a_1,\cdots,a_{q-1})=(-1)^{q-1}\frac{\partial}{\partial\la}\g_{\la_1,\cdots,\la_{q-1},\la}(a_1,\cdots,a_{q-1},L)|_{\la=0}.
\end{eqnarray}
By direct computations (referring to \cite{YW}), we have
\begin{align}\label{al3.6}
((d\tau+&\tau d)\g)_{\lambda_1,\cdots,\lambda_{q}}(a_1,\cdots,a_q)\nonumber\\
=&(-1)^q\frac{\partial}{\partial\lambda}\sum_{i=1}^q(-1)^{i+q+1}\g_{\lambda_i+\lambda,\lambda_1,\cdots,\hat{\lambda_i},\cdots,\lambda_{q}}([{a_i}_{\lambda_i}L],a_1,\cdots,\hat{a_i},\cdots,a_{q})|_{\lambda=0}\nonumber\\
=&\frac{\partial}{\partial\lambda}\sum_{i=1}^q\g_{\lambda_1,\cdots,\lambda_{i-1},\lambda_i+\lambda,\lambda_{i+1},\cdots,\lambda_{q}}(a_1,\cdots,a_{i-1},[{a_i}_{\lambda_i}L],a_{i+1},\cdots,a_{q})|_{\lambda=0}\nonumber\\
=&\frac{\partial}{\partial\lambda}\sum_{i=1}^k(\lambda_i-\lambda)\g_{\lambda_1,\cdots,\lambda_{i-1},\lambda_i+\lambda,\lambda_{i+1},\cdots,\lambda_{q}}(a_1,\cdots,a_{i-1},a_i,a_{i+1},\cdots,a_{q})|_{\lambda=0}\nonumber\\
&+\frac{\partial}{\partial\lambda}\sum_{i=k+1}^{k+l}(\lambda_i-\frac{1}{2}\lambda)\g_{\lambda_1,\cdots,\lambda_{i-1},\lambda_i+\lambda,\lambda_{i+1},\cdots,\lambda_{q}}(a_1,\cdots,a_{i-1},a_i,a_{i+1},\cdots,a_{q})|_{\lambda=0}\nonumber\\
&+\frac{\partial}{\partial\lambda}\sum_{i=k+l+1}^{k+l+m}\lambda_i\g_{\lambda_1,\cdots,\lambda_{i-1},\lambda_i+\lambda,\lambda_{i+1},\cdots,\lambda_{q}}(a_1,\cdots,a_{i-1},a_i,a_{i+1},\cdots,a_{q})|_{\lambda=0}\nonumber\\
=&({\rm deg\,} \g-k-\frac{l}{2})\g,
\end{align}
where ${\rm deg}\, \g$ is the total degree of $\g$ in $\la_1,\cdots,\la_{q}$.
Therefore, if a $q$-cocyle $\g$ satisfies ${\rm deg}\, \g\neq k+\frac{l}{2}$, it must be a coboundary.  Only those homogeneous cochains whose degree as a polynomial is equal to $k+\frac{l}{2}$ contribute to the cohomology of $\t C^\bullet(\SV,\C)$.

Consider the quadratic inequality
\begin{equation}\label{ineq}
\begin{pmatrix}k\\2\end{pmatrix}+\begin{pmatrix}l\\2\end{pmatrix}+\begin{pmatrix}m\\2\end{pmatrix}\leq k+\frac{l}{2}.
\end{equation}

\begin{lemm}\label{le0}
All non-negative integral solutions satisfying $k+\frac{l}{2}\in\Z_+$ of the inequality (\ref{ineq}) are listed in the following table:
\end{lemm}
%%%%%%%%%%%%%%
\begin{CJK*}{GBK}{song}
\setlength{\LTleft}{100pt} \setlength{\LTright}{100pt} %表格与页面左右边缘之间的矩离均为０
\begin{longtable}{|c|c|c|c|c|c|}\hline
			$q=k+l+m$&$(k,l,m)$&$\begin{pmatrix}k\\2\end{pmatrix}+\begin{pmatrix}l\\2\end{pmatrix}+\begin{pmatrix}m\\2\end{pmatrix}$&
${\rm deg}\, \g=k+\frac{l}{2}$\\\hline
			$0$&$(0,0,0)$&$0$&$0$\\\hline
                    $1$&$(0,0,1)$&$0$&$0$\\
                     &$(1,0,0)$&$0$&$1$\\\hline
                     &$(0,2,0)$&$1$&$1$\\
                     $2$&$(1,0,1)$&$0$&$1$\\
                    &$(2,0,0)$&$1$&$2$\\\hline
                     &$(0,2,1)$&$1$&$1$\\
                    &$(1,0,2)$&$1$&$1$\\
                     $3$&$(1,2,0)$&$1$&$2$\\
                    &$(2,0,1)$&$1$&$2$\\
                     &$(3,0,0)$&$3$&$3$\\\hline
                     &$(1,2,1)$&$1$&$2$\\
                    $4$&$(2,0,2)$&$2$&$2$\\
                     &$(2,2,0)$&$2$&$3$\\
                     &$(3,0,1)$&$3$&$3$\\\hline
                     &$(1,2,2)$&$2$&$2$\\
                     $5$&$(2,2,1)$&$2$&$3$\\
                      &$(3,2,0)$&$4$&$4$\\\hline
                    $6$&$(2,2,2)$&$3$&$3$\\
                     &$(3,2,1)$&$4$&$4$\\\hline
	\end{longtable}
\end{CJK*}
	
By the conclusion of the above table, we can obtain the following result immediately.
\begin{lemm}
$\t H^q(\SV,\C)=0, q\geq7$.
\end{lemm}
\begin{lemm}
$\t H^1(\SV,\C)=0$.
\end{lemm}
{\it Proof.} By the Lemma \ref{le0}, we only need to consider $(k,l,m)=(1,0,0)$ and $(0,0,1)$. Let $\gamma$ be a $1$-cocyle. Assume $\gamma_\lambda(L)=a\lambda$,  $\gamma_\lambda(Y)=0$ and $\gamma_\lambda(M)=b$. Then $d\gamma_{\lambda_1,\lambda_2}(L,L)=-a(\lambda_1-\lambda_2)(\lambda_1+\lambda_2)=0$ and $d\gamma_{\lambda_1,\lambda_2}(L,M)=b\lambda_2=0$, implying $\gamma=0$.
\QED

\begin{lemm}
$\t H^2(\SV,\C)=0$.
\end{lemm}
{\it Proof.} We only need to consider $(k,l,m)=(0,2,0), (1,0,1)$ and $(2,0,0)$. Let $\gamma$ be a $2$-cocyle. Assume $\gamma_{\lambda_1,\lambda_2}(L, L)=a(\lambda_1+\lambda_2)(\lambda_1-\lambda_2)$, $\gamma_{\lambda_1,\lambda_2}(L, M)=b\lambda_1+c\lambda_2$ and $\gamma_{\lambda_1,\lambda_2}(Y,Y)=e(\lambda_1-\lambda_2)$ for some $a,b,c, e\in\C$.
Let $\phi$ be a $1$-cochain defined by $\phi_\lambda(L)=a\lambda, \phi_\lambda(Y)=0, \phi_\lambda(M)=e$.
Replacing $\gamma$ by $\gamma+d\phi$, we can assume that $\gamma_{\lambda_1,\lambda_2}(L, L)=0$ and $\gamma_{\lambda_1,\lambda_2}(Y,Y)=0$.
Then by  $$d\gamma_{\lambda_1,\lambda_2, \lambda_3}(L,Y,Y)=-\gamma_{\lambda_2+ \lambda_3,\lambda_1}((\partial+2\lambda_2)M,L)=
(\lambda_2-\lambda_3)(b\lambda_1+c\lambda_2+c\lambda_3)=0,$$ we have $\gamma=0$.
\QED

\begin{lemm}$\t H^3(\SV,\C)=\C \Phi$, where $\Phi$ is defined by $$\Phi_{\lambda_1,\lambda_2,\lambda_3}(L, L, L)=(\lambda_1-\lambda_2)(\lambda_1-\lambda_3)(\lambda_2-\lambda_3).$$ In particular,
${\rm dim}\,\t H^3(\SV,\C)=1$.
\end{lemm}
{\it Proof.} We should consider $(k,l,m)=(3,0,0), (2,0,1), (1,2,0), (1,0,2)$ and $(0,2,1)$.
Let $\gamma$ be an arbitrary $3$-cocycle.
Then we can assume
\begin{eqnarray*}
&&\gamma_{\lambda_1,\lambda_2,\lambda_3}(L, L, L)=a(\lambda_1-\lambda_2)(\lambda_1-\lambda_3)(\lambda_2-\lambda_3),\\&&
\gamma_{\lambda_1,\lambda_2,\lambda_3}(L, L, M)=(\lambda_1-\lambda_2)(b\lambda_1+b\lambda_2+c\lambda_3),\\
&&\gamma_{\lambda_1,\lambda_2,\lambda_3}(L,Y,Y)=
(\lambda_2-\lambda_3)(e\lambda_1+f\lambda_2+f\lambda_3),\\
&&\gamma_{\lambda_1,\lambda_2,\lambda_3}(L,M,M)=
g(\lambda_2-\lambda_3),\\
&&\gamma_{\lambda_1,\lambda_2,\lambda_3}(Y,Y,M)=
h(\lambda_1-\lambda_2),
\end{eqnarray*}
where $a,b,c,e,f,g,h\in\C$. In the following, we try to determine these parameters.

%Let $\phi$ be a $2$-cochain defined by $\phi_{\lambda_1,\lambda_2}(M,M)=1$.
%Then
%\begin{eqnarray*}
%&&(d\phi)_{\lambda_1,\lambda_2,\lambda_3}(L,M,M)=\lambda_2-\lambda_3,\\
%&&(d\phi)_{\lambda_1,\lambda_2,\lambda_3}(Y,Y,M)=-(\lambda_1-\lambda_2).
%\end{eqnarray*} Replacing $\gamma$ by $\gamma+h(d\phi)$, we can assume $h=0$ and $\gamma_{\lambda_1,\lambda_2,\lambda_3}(Y,Y,M)=0$.
%Then
%\begin{eqnarray*}
%&(d\gamma)_{\lambda_1,\lambda_2,\lambda_3,\lambda_4}(L,Y,Y,M)\\=&
%-\gamma_{\lambda_2+\lambda_3,\lambda_1,\lambda_4}((\partial+2\lambda_2)M,L,M)\\
%=&g(\lambda_2-\lambda_3)(\lambda_2+\lambda_3-\lambda_4)=0,
%\end{eqnarray*}
%which implies $g=0$.

Let $\varphi$ be a $2$-cochain defined by $\varphi_{\lambda_1,\lambda_2}(L,M)=b\lambda_1$ and $\varphi_{\lambda_1,\lambda_2}(Y,Y)=-f(\lambda_1-\lambda_2)$.
Then
\begin{eqnarray*}
&&(d\varphi)_{\lambda_1,\lambda_2,\lambda_3}(L,L,M)=
(\lambda_1-\lambda_2)(-b\lambda_1-b\lambda_2-b\lambda_3),\\
&&(d\varphi)_{\lambda_1,\lambda_2,\lambda_3}(L,Y,Y)=
(\lambda_2-\lambda_3)(b\lambda_1-f\lambda_2-f\lambda_3).
\end{eqnarray*}
Replacing $\gamma$ by $\gamma+d\varphi$, we can assume $b=f=0$.  Consequently, $\gamma_{\lambda_1,\lambda_2,\lambda_3}(L,L,M)=
c(\lambda_1-\lambda_2)\lambda_3$ and $\gamma_{\lambda_1,\lambda_2,\lambda_3}(L,Y,Y)=
e(\lambda_2-\lambda_3)\lambda_1$ for some $c,e\in \C$.
Then by direct computations, we have
\begin{eqnarray*}
(d\gamma)_{\lambda_1,\lambda_2,\lambda_3,\lambda_4}(L,L,Y,Y)
=-(\lambda_1-\lambda_2)(\lambda_3-\lambda_4)
(e\lambda_1+e\lambda_2+(c+e)\lambda_3+(c+e)\lambda_4)=0,
\end{eqnarray*}
which implies $c=e=0$.

%Let $\psi$ be a $2$-cochain defined by $\varphi_{\lambda_1,\lambda_2}(Y,Y)=\lambda_1-\lambda_2$.
%Then
%\begin{eqnarray*}
%(d\varphi)_{\lambda_1,\lambda_2,\lambda_3}(L,Y,Y)=
%\lambda_1(\lambda_2-\lambda_3).
%\end{eqnarray*}

Similarly, we can obtain that
\begin{eqnarray*}
(d\gamma)_{\lambda_1,\lambda_2,\lambda_3,\lambda_4}(L,Y,Y,M)
=(\lambda_2-\lambda_3)((g+h)\lambda_2+(g+h)\lambda_3+(h-g)\lambda_4)=0,
\end{eqnarray*}
which implies $g=h=0$.

By the cohomology theory of the Virasoro conformal algebra in \cite{CK}, we know that the cochain defined by $\Phi_{\lambda_1,\lambda_2,\lambda_3}(L, L, L)=(\lambda_1-\lambda_2)(\lambda_1-\lambda_3)(\lambda_2-\lambda_3)$ is a cocycle, but not a coboundary. Thus we have $\t H^3(\SV,\C)=\C \Phi$.
\QED

\begin{rema}
The above skew-symmetric function $\Phi: \SV\otimes\SV\otimes\SV\to\C[\lambda_1,\lambda_2,\lambda_3]$ has values $(\lambda_1-\lambda_2)(\lambda_1-\lambda_3)(\lambda_2-\lambda_3)$ on $L\otimes L\otimes L$ and $0$ on others.
 \end{rema}
%Let $\gamma$ be a $3$-cochain defined by  %$\gamma_{\lambda_1,\lambda_2,\lambda_3}(L, M, M)=\lambda_2-\lambda_3$. It %was shown in \cite{BKV} that $\gamma$ is a $3$-cocycle and not a %$3$-coboundary.

\begin{lemm}$\t H^4(\SV,\C)=0$.
\end{lemm}
{\it Proof.} For $q=4$, we should consider $(k,l,m)=(3,0,1), (2,2,0), (2,0,2)$ and $(1,2,1)$.
Let $\gamma$ be an arbitrary $4$-cocycle.
Then we can assume $\gamma$ is defined by
\begin{eqnarray*}
&&\gamma_{\lambda_1,\lambda_2,\lambda_3,\lambda_4}(L, L, L,M)=a(\lambda_1-\lambda_2)(\lambda_1-\lambda_3)(\lambda_2-\lambda_3),\\&&
\gamma_{\lambda_1,\lambda_2,\lambda_3,\lambda_4}(L, L, M,M)=b(\lambda_1-\lambda_2)(\lambda_3-\lambda_4),\\
&&\gamma_{\lambda_1,\lambda_2,\lambda_3,\lambda_4}(L,L,Y,Y)=
(\lambda_1-\lambda_2)(\lambda_3-\lambda_4)
(c\lambda_1+c\lambda_2+e\lambda_3+e\lambda_4),\\
&&\gamma_{\lambda_1,\lambda_2,\lambda_3,\lambda_4}(L,Y,Y,M)=
(\lambda_2-\lambda_3)(f\lambda_1+g\lambda_2+g\lambda_3+h\lambda_4),
\end{eqnarray*}
where $a,b,c,e,f,g,h\in\C$.

Let $\phi$ be a $3$-cochain defined by
$\phi_{\lambda_1,\lambda_2,\lambda_3}(L,L,M)=
(\lambda_1-\lambda_2)\lambda_3$.
Then by direct computations, we have
\begin{eqnarray*}
&&(d\phi)_{\lambda_1,\lambda_2,\lambda_3,\lambda_4}(L,L,Y,Y)
=-(\lambda_1-\lambda_2)(\lambda_3-\lambda_4)
(\lambda_3+\lambda_4),\\
&&(d\phi)_{\lambda_1,\lambda_2,\lambda_3,\lambda_4}(L,L,L,M)
=0.
\end{eqnarray*}
Let $\psi$ be a $3$-cochain defined by
$\psi_{\lambda_1,\lambda_2,\lambda_3}(L,Y,Y)=
\lambda_1(\lambda_2-\lambda_3)$.
Then by direct computations, we have
\begin{eqnarray*}
(d\psi)_{\lambda_1,\lambda_2,\lambda_3,\lambda_4}(L,L,Y,Y)
=-(\lambda_1-\lambda_2)(\lambda_3-\lambda_4)
(\lambda_1+\lambda_2+\lambda_3+\lambda_4).
\end{eqnarray*}
Replacing $\gamma$ by $\gamma+(e-c)(d\phi)+c(d\psi)$, we can assume $c=e=0$ and $\gamma_{\lambda_1,\lambda_2,\lambda_3,\lambda_4}(L,L,Y,Y)=0$.
Then by
\begin{eqnarray*}
(d\gamma)_{\lambda_1,\lambda_2,\lambda_3,\lambda_4,\lambda_5}(L,L,L,Y,Y)
=a(\lambda_1-\lambda_2)(\lambda_1-\lambda_3)(\lambda_2-\lambda_3)(\lambda_4-\lambda_5)=0,
\end{eqnarray*}
we deduce that $a=0$. Therefore, $\gamma_{\lambda_1,\lambda_2,\lambda_3,\lambda_4}(L,L,L,M)=0$.

Similarly, let $\Phi$ be a $3$-cochain defined by
$\Phi_{\lambda_1,\lambda_2,\lambda_3}(L,M,M)=
\lambda_2-\lambda_3$.
Then by direct computations, we have
\begin{eqnarray*}
&&(d\Phi)_{\lambda_1,\lambda_2,\lambda_3,\lambda_4}(L,Y,Y,M)
=(\lambda_2-\lambda_3)(\lambda_2+\lambda_3-\lambda_4),\\
&&(d\Phi)_{\lambda_1,\lambda_2,\lambda_3,\lambda_4}(L,L,M,M)
=0.
\end{eqnarray*}
Let $\Psi$ be a $3$-cochain defined by
$\Psi_{\lambda_1,\lambda_2,\lambda_3}(Y,Y,M)=
\lambda_1-\lambda_2$.
Then by direct computations, we have
\begin{eqnarray*}
(d\Psi)_{\lambda_1,\lambda_2,\lambda_3,\lambda_4}(L,Y,Y,M)
=(\lambda_2-\lambda_3)(\lambda_2+\lambda_3+\lambda_4).
\end{eqnarray*}
Replacing $\gamma$ by $\gamma-\frac{g-h}{2}(d\Phi)-\frac{g+h}{2}(d\Psi)$, we can assume $g=h=0$.

Now, we can assume
\begin{eqnarray*}
&&\gamma_{\lambda_1,\lambda_2,\lambda_3,\lambda_4}(L, L, M,M)=b(\lambda_1-\lambda_2)(\lambda_3-\lambda_4),\\
&&\gamma_{\lambda_1,\lambda_2,\lambda_3,\lambda_4}(L,Y,Y,M)=
f\lambda_1(\lambda_2-\lambda_3),
\end{eqnarray*}
where $b,f\in\C$.
Then by \begin{eqnarray*}
&&(d\gamma)_{\lambda_1,\lambda_2,\lambda_3,\lambda_4,\lambda_5}(L,L,Y,Y,M)\\
=&&b(\lambda_1-\lambda_2)(\lambda_3-\lambda_4)(-\lambda_3-\lambda_4+\lambda_5)
-
f(\lambda_1-\lambda_2)(\lambda_3-\lambda_4)
(\lambda_1+\lambda_2+\lambda_3+\lambda_4+\lambda_5),
\end{eqnarray*}
we deduce that $b=f=0$. Thus we have $\gamma=0$.\QED

\begin{lemm}\label{l7}
$\t H^5(\SV,\C)=\C \Lambda\oplus \C\Psi$, where $$\Lambda_{\lambda_1,\lambda_2,\lambda_3,\lambda_4,\lambda_5}(L, L, Y,Y,M)=(\lambda_1-\lambda_2)(\lambda_3-\lambda_4)\lambda_5$$ and $$\Psi_{\lambda_1,\lambda_2,\lambda_3,\lambda_4,\lambda_5}(L,Y,Y,M,M)=
(\lambda_2-\lambda_3)(\lambda_4-\lambda_5).$$ In particular,
${\rm dim}\,\t H^5(\SV,\C)=2$.
\end{lemm}
{\it Proof.} For $q=5$, we should consider $(k,l,m)=(3,2,0), (2,2,1)$ and $(1,2,2)$. Let $\gamma$ be an arbitrary $5$-cocycle.
Then we can assume that $\gamma$ is defined by
\begin{eqnarray*}
&&\gamma_{\lambda_1,\lambda_2,\lambda_3,\lambda_4,\lambda_5}(L, L, L,Y,Y)=a(\lambda_1-\lambda_2)(\lambda_1-\lambda_3)
(\lambda_2-\lambda_3)(\lambda_4-\lambda_5),\\&&
\gamma_{\lambda_1,\lambda_2,\lambda_3,\lambda_4,\lambda_5}(L, L, Y,Y,M)=(\lambda_1-\lambda_2)(\lambda_3-\lambda_4)
(b\lambda_1+b\lambda_2+c\lambda_3+c\lambda_4+e\lambda_5),\\
&&\gamma_{\lambda_1,\lambda_2,\lambda_3,\lambda_4, \lambda_5}(L,Y,Y,M,M)=f
(\lambda_2-\lambda_3)(\lambda_4-\lambda_5),
\end{eqnarray*}
where $a,b,c,e,f\in\C$.

Let $\phi$ be a $4$-cochain defined by
$$\phi_{\lambda_1,\lambda_2,\lambda_3,\lambda_4}(L, L, L,M)=(\lambda_1-\lambda_2)(\lambda_1-\lambda_3)(\lambda_2-\lambda_3).$$
Then $d\phi$ is given by \begin{eqnarray*}
&&(d\phi)_{\lambda_1,\lambda_2,\lambda_3,\lambda_4,\lambda_5}(L,L,L,Y,Y)
=(\lambda_1-\lambda_2)(\lambda_1-\lambda_3)
(\lambda_2-\lambda_3)(\lambda_4-\lambda_5),\\
&&(d\phi)_{\lambda_1,\lambda_2,\lambda_3,\lambda_4,\lambda_5}(L,L,L,L,M)=0.
\end{eqnarray*}
Replacing $\gamma$ by $\gamma-a(d\phi)$, we have
$\gamma_{\lambda_1,\lambda_2,\lambda_3,\lambda_4,\lambda_5}(L,L,L,Y,Y)=0$.

Define
\begin{eqnarray*}
&&\psi_{\lambda_1,\lambda_2,\lambda_3,\lambda_4}(L, L, M,M)=(c-b)(\lambda_1-\lambda_2)(\lambda_3-\lambda_4),\\
&&\psi_{\lambda_1,\lambda_2,\lambda_3,\lambda_4}(L,Y,Y,M)=
b\lambda_1(\lambda_2-\lambda_3).
\end{eqnarray*}
Then \begin{eqnarray*}
(d\psi)_{\lambda_1,\lambda_2,\lambda_3,\lambda_4,\lambda_5}(L,L,Y,Y,M)
=-(\lambda_1-\lambda_2)(\lambda_3-\lambda_4)
(b\lambda_1+b\lambda_2+c\lambda_3+c\lambda_4+(2b-c)\lambda_5).
\end{eqnarray*}
Replacing $\gamma$ by $\gamma-d(\psi)$, we can assume $$\gamma_{\lambda_1,\lambda_2,\lambda_3,\lambda_4,\lambda_5}(L, L, Y,Y,M)=e(\lambda_1-\lambda_2)(\lambda_3-\lambda_4)
\lambda_5.$$
Then by
\begin{eqnarray*}
(d\gamma)_{\lambda_1,\lambda_2,\lambda_3,\lambda_4,\lambda_5,\lambda_6}(L,L,L,Y,Y,M)
=0.
\end{eqnarray*}
Define $\Lambda_{\lambda_1,\lambda_2,\lambda_3,\lambda_4,\lambda_5}(L, L, Y,Y,M)=(\lambda_1-\lambda_2)(\lambda_3-\lambda_4)
\lambda_5$, which is a cocycle, but not a coboundary.

Similarly, we can show that a cochain defined by $\Psi_{\lambda_1,\lambda_2,\lambda_3,\lambda_4}(L,Y,Y,M,M)=
(\lambda_2-\lambda_3)(\lambda_4-\lambda_5)$ is a cocycle, but not a coboundary.
\QED

\begin{lemm}$\t H^6(\SV,\C)=\C \Omega\oplus \C\Theta$, where $$\Omega_{\lambda_1,\lambda_2,\lambda_3,\lambda_4,\lambda_5,\lambda_6}(L, L, L,Y,Y,M)=(\lambda_1-\lambda_2)(\lambda_1-\lambda_3)
(\lambda_2-\lambda_3)(\lambda_4-\lambda_5)$$ and $$\Theta_{\lambda_1,\lambda_2,\lambda_3,\lambda_4,\lambda_5,\lambda_6}(L, L, Y,Y,M,M)=(\lambda_1-\lambda_2)
(\lambda_3-\lambda_4)(\lambda_5-\lambda_6).$$ In particular,
${\rm dim}\,\t H^6(\SV,\C)=2$.
\end{lemm}
{\it Proof.}For $q=6$, we should consider $(k,l,m)=(3,2,1)$ and $(2,2,2)$.

First, let $\gamma$ be the $6$-cochain defined by \begin{eqnarray*}
&&\gamma_{\lambda_1,\lambda_2,\lambda_3,\lambda_4,\lambda_5,\lambda_6}(L, L, L,Y,Y,M)=(\lambda_1-\lambda_2)(\lambda_1-\lambda_3)
(\lambda_2-\lambda_3)(\lambda_4-\lambda_5).
\end{eqnarray*} It can be only a coboundary of $\phi$,  defined by
\begin{align*}
\phi_{\lambda_1,\lambda_2,\lambda_3,\lambda_4,\lambda_5}(L,L,Y,Y,M)\quad and\quad  \phi_{\lambda_1,\lambda_2,\lambda_3,\lambda_4,\lambda_5}(L,L,L,M,M).
\end{align*}
But by the proof of Lemma \ref{l7}, $\phi_{\lambda_1,\lambda_2,\lambda_3,\lambda_4,\lambda_5}(L,L,Y,Y,M)$ is a cocycle. Furthermore, the degree of $\phi_{\lambda_1,\lambda_2,\lambda_3,\lambda_4,\lambda_5}(L,L,L,M,M)$ is at least $4$, which implies $\gamma$ is not a coboundary. It is straightforward to check $d\gamma=0$.

Second, let $\gamma$ be the $6$-cochain defined by \begin{eqnarray*}
&&\gamma_{\lambda_1,\lambda_2,\lambda_3,\lambda_4,\lambda_5,\lambda_6}(L, L, Y,Y,M,M)=(\lambda_1-\lambda_2)
(\lambda_3-\lambda_4)(\lambda_5-\lambda_6).
\end{eqnarray*} It can be only a coboundary of $\phi$,  defined by
\begin{align*}
\phi_{\lambda_1,\lambda_2,\lambda_3,\lambda_4,\lambda_5}(L,Y,Y,M,M)\quad and \quad \phi_{\lambda_1,\lambda_2,\lambda_3,\lambda_4,\lambda_5}(L,L,M,M,M).
\end{align*}
As shown above, we can obtain that $\phi_{\lambda_1,\lambda_2,\lambda_3,\lambda_4,\lambda_5}(L,L,Y,Y,M)$ is a cocycle and the degree of $\phi_{\lambda_1,\lambda_2,\lambda_3,\lambda_4,\lambda_5}(L,L,M,M,M)$ is at least $4$. Thus $\gamma$ is not a coboundary. Moreover, the degree of a $7$-cochain  is at least $5$ on $(L,L,L,Y,Y,M,M)$ and at least $7$ on $(L,L,Y,Y,Y,Y,M)$. Then we deduce that $d\gamma=0$.
\QED

 As mentioned above, we can obtain the main result of this section as follows.
\begin{theo}\label{th1}The dimension of ${\t H}^q (\mathcal{SV},\C)$ is given by
\begin{eqnarray*}
{\rm dim\,\t H}^q(\SV,\C)=\left\{
\begin{array}{ll}
1 &{\mbox if}\ q=0,3,\\
2 &{\mbox if}\ q=5,6,\\
0 &{\mbox otherwise}.
\end{array}
\right.
\end{eqnarray*}
In particular,
\begin{eqnarray}\label{hd001}
{\t H^q(\mathcal{SV}, \C)}=\left\{
\begin{array}{ll}
\C  &{\mbox if}\ q=0,\\
\C (\Phi) &{\mbox if}\ q=3,\\
\C (\Lambda)\oplus \C(\Psi) &{\mbox if}\ q=5,\\
\C (\Omega)\oplus \C(\Theta) &{\mbox if}\ q=6,\\
0 &{\mbox otherwise},
\end{array}
\right.
\end{eqnarray}
where
\begin{align*}
\Phi&_{\lambda_1,\lambda_2,\lambda_3}(L, L, L)=(\lambda_1-\lambda_2)(\lambda_1-\lambda_3)(\lambda_2-\lambda_3),\\
\Lambda&_{\lambda_1,\lambda_2,\lambda_3,\lambda_4,\lambda_5}(L,L,Y,Y,M)=(\lambda_1-\lambda_2)(\lambda_3-\lambda_4)\lambda_5,\\
\Psi&_{\lambda_1,\lambda_2,\lambda_3,\lambda_4,\lambda_5}(L,Y,Y,M,M)=(\lambda_2-\lambda_3)(\lambda_4-\lambda_5),\\
\Omega&_{\lambda_1,\lambda_2,\lambda_3,\lambda_4,\lambda_5,\lambda_6}(L, L, L,Y,Y,M)=(\lambda_1-\lambda_2)(\lambda_1-\lambda_3)
(\lambda_2-\lambda_3)(\lambda_4-\lambda_5),\\
\Theta&_{\lambda_1,\lambda_2,\lambda_3,\lambda_4,\lambda_5,\lambda_6}(L, L, Y,Y,M,M)=(\lambda_1-\lambda_2)(\lambda_3-\lambda_4)(\lambda_5-\lambda_6).
\end{align*}
\end{theo}
{\it Proof.} It follows directly from the previous Lemmas.
\QED
\begin{rema}
The corresponding annihilation algebra of $\mathcal{SV}$ is
\begin{eqnarray*}
&&\textit{Lie}(\mathcal{SV})^+= \sum_{m\geq -1}\CC L_m+\sum_{n\geq 0}\CC M_n +\sum_{p\in\frac{1}{2}+\Z_+}\CC Y_{p},\\
&&[L_m,L_{n}]=(m-n)L_{m+n},\ \ \
[L_m,M_n]=-nM_{m+n},\\
&&[\,Y_p\,,Y_{q}\,]=(p-q)M_{p+q},\ \ \ \ \
[\,L_m,Y_p\,]=(\frac{m}{2}-p)Y_{m+p},
\end{eqnarray*}
which is a 'half part' of the Schr\"odinger-Virasoro Lie algebra (\cite{SY}). So by \cite{BKV}, the dimension of all the cohomology groups of $\textit{Lie}(\mathcal{SV})^+$ are given by
\begin{eqnarray*}
{\rm dim\, H}^q(\textit{Lie}(\mathcal{SV})^+,\C)=\left\{
\begin{array}{ll}
1 &{\mbox if}\ q=0,3,\\
2 &{\mbox if}\ q=5,6,\\
0 &{\mbox otherwise}.
\end{array}
\right.
\end{eqnarray*}
 \end{rema}

\section{Reduced cohomology of $\mathcal{SV}$ with coefficients in its module}
In this section, we compute the reduced cohomology groups of $\mathcal{SV}$ with coefficients in its trivial module $\mathbb{C}_a$ and with coefficients in $M_{\alpha,\beta}$ in the case of $\beta\neq0$.

%The group $H^q (\mathcal{SV},\C)$ varies depending on the values of $a$. So we divide into two parts.
\subsection{Computation of $H^q (\mathcal{SV},\C)$}
\begin{theo}\label{th2} The dimension of $H^q (\mathcal{SV},\C)$ are given by
\begin{eqnarray*}
{\rm dim\, H}^q(\SV,\C)=\left\{
\begin{array}{ll}
1 &{\mbox if}\ q=0,2,3,\\
2 &{\mbox if}\ q=4,6,\\
4 &{\mbox if}\ q=5,\\
0 &{\mbox otherwise}.
\end{array}
\right.
\end{eqnarray*}
\end{theo}
{\it Proof.}
By Proposition \ref{pro2.12}, the map $\g \mapsto \partial \cdot \g$ gives an isomorphism such that $\t H^q(\mathcal{SV}, \C)\cong H^q(\partial\t C^\bullet) $ for all $q\geq1$. Therefore, we can obtain the following result immediately by the discussion of Section 3.1.

\begin{eqnarray}\label{hd}
{H^q(\partial\t C^\bullet) }=\left\{
\begin{array}{llllll}
\C (\partial \Phi) &{\mbox if}\ q=3,\\
\C (\partial \Lambda)\oplus \C(\partial \Psi) &{\mbox if}\ q=5,\\
\C (\partial \Omega)\oplus \C(\partial \Theta) &{\mbox if}\ q=6,\\
0 &{\mbox otherwise}.
\end{array}
\right.
\end{eqnarray}

Similar to the discussions in Section 2, we can obtain the following long exact sequence of cohomology groups:
\begin{align}\label{longSVC}
\cdots\longrightarrow& \ \ H^q(\partial\t C^\bullet) \   \stackrel{i_q}{\longrightarrow} \ \ \t H^q(\mathcal{SV},\C)\  \stackrel{p_q}{\longrightarrow} \ \ H^q (\mathcal{SV},\C) \ \stackrel{w_q}\longrightarrow\\
\longrightarrow& H^{q+1}(\partial\t C^\bullet) \stackrel{i_{q+1}}{\longrightarrow} \t H^{q+1}(\mathcal{SV},\C) \stackrel{p_{q+1}}{\longrightarrow} H^{q+1} (\mathcal{SV},\C) \longrightarrow\cdots\nonumber
\end{align}
where $i_q, p_q$ are induced by $i,p$ in (\ref{longexact}) respectively and $w_q$ is the $q$-th connecting homomorphism. Take $\partial \g \in H^q(\partial\t C^\bullet) $ with a nonzero element $\g\in\t H^q(\mathcal{SV},\C)$ of degree $k+\frac{l}{2}$, we can obtain that $i_q(\partial \g)=\partial \g\in \t H^q(\mathcal{SV},\C)$.
Since $\partial\gamma=(\sum\lambda_i)\gamma$, we have $${\rm deg}\,(\partial \g)={\rm deg}\,(\g)+1=k+\frac{l}{2}+1.$$
Thus, $\partial \g=0 \in\t H^q(\mathcal{SV},\C)$ (Here we use the property $a=0$).%  because
%$${\rm deg}\,(\partial \g)={\rm deg}\,(\g)+1=k+\frac{l}{2}+1.$$
Thus, the image of $i_q$ is zero. %, i.e., $im(i_q)=0$ for any $q\in\mathbb{Z_+}$. Since $ker(p_q)=im(i_q)=0$ and $im(w_q)=ker(i_{q+1})= H^{q+1}(\partial\t C^\bullet)$, we can get
Then the long exact sequence (\ref{longSVC}) splits into
the following short exact sequence immediately:
\begin{align}\label{iq}
0\longrightarrow \t H^q(\mathcal{SV},\C)\stackrel{p_q}{\longrightarrow} H^q (\mathcal{SV},\C)\stackrel{w_q}\longrightarrow H^{q+1}(\partial\t C^\bullet)\longrightarrow 0.
\end{align}
for all $q\geq1$. Thus, we have
\begin{align}\label{dimc}
{\rm dim}\, H^q (\mathcal{SV},\C) =& {\rm dim}\, \t H^q(\mathcal{SV},\C)+ {\rm dim}\, H^{q+1}(\partial\t C^\bullet)\nonumber \\
=&{\rm dim}\,\t H^q(\mathcal{SV},\C)+{\rm dim}\, \t H^{q+1}(\mathcal{SV},\C)
\end{align}
for all $q\geq1$. Then the result follows.\QED

\begin{rema}It was shown in \cite{SY} that there is a unique nontrivial extension of $\SV$ by a $1$-dimensional center.
This coincides with our result ${\rm dim\, H}^2(\SV,\C)=1$.
\end{rema}

It is not difficult to check by (\ref{iq}) that the basis of $H^q(\SV,\C)$ can be obtained by combining the images of a basis of $\t H^q(\mathcal{SV},\C)$ with the pre-images of a basis of $\t H^{q+1}(\mathcal{SV},\C)$. Let $\g$ be a nonzero $(q+1)$-cocycle of degree $k+\frac{l}{2}$ such that $\partial\g \in H^{q+1}(\partial\t C^\bullet)$. By (\ref{al3.6}), we can obtain that
\begin{align}\label{dimcc}
d(\tau(\partial\g))=(d\tau+\tau d)(\partial\g)=(deg(\partial\g)-k-\frac{l}{2})(\partial\g)=((k+\frac{l}{2}+1)-k-\frac{l}{2})(\partial\g)=\partial\g.
\end{align}
Thus, the pre-image $w_q^{-1}(\partial\g)$ of $\partial\g$ under the connecting homomorphism $w_q$ is $\tau(\partial\g)$, i.e., $w_q^{-1}(\partial\g)=\tau(\partial\g).$

Due to the above discussions, we can obtain a specific system to calculate the basis of $H^q(\SV,\C)$ so that we can describe the structure of $H^q(\SV,\C)$ more clearly.

\begin{coro}
\begin{eqnarray}\label{hd2}
{ H^q(\mathcal{SV}, \C)}=\left\{
\begin{array}{ll}
\C  &{\mbox if}\ q=0,\\
\C (\bar{\Phi}) &{\mbox if}\ q=2,\\
\C (\Phi) &{\mbox if}\ q=3,\\
\C (\bar{\Lambda})\oplus \C(\bar{\Psi}) &{\mbox if}\ q=4,\\
\C (\Lambda)\oplus \C(\Psi)\oplus\C (\bar{\Omega})\oplus \C(\bar{\Theta}) &{\mbox if}\ q=5,\\
\C (\Omega)\oplus \C(\Theta) &{\mbox if}\ q=6,\\
0 &{\mbox otherwise},
\end{array}
\right.
\end{eqnarray}
where $\bar{X}=\tau(\partial X)$ and $X\in \{\Phi,\Lambda,\Psi,\Omega,\Theta\}$ as shown in Theorem \ref{th1}. More specifically,
\begin{align*}
\bar{\Phi}&_{\lambda_1,\lambda_2}(L, L)=-\lambda_1^3+\lambda_2^3,\\
\bar{\Lambda}&_{\lambda_1,\lambda_2,\lambda_3,\lambda_4}(L,Y,Y,M)=(\lambda_2+\lambda_3+\lambda_4)(\lambda_2-\lambda_3)\lambda_4,\\
\bar{\Psi}&_{\lambda_1,\lambda_2,\lambda_3,\lambda_4}(Y,Y,M,M)=(\lambda_1-\lambda_2)(\lambda_3-\lambda_4),\\
\bar{\Omega}&_{\lambda_1,\lambda_2,\lambda_3,\lambda_4,\lambda_5}(L,L,Y,Y,M)=(\lambda_1-\lambda_2)(\lambda_3-\lambda_4)
(\lambda_1\lambda_2-(\lambda_1+\lambda_2)(\lambda_1+\lambda_2+\lambda_3+\lambda_4+\lambda_5)),\\
\bar{\Theta}&_{\lambda_1,\lambda_2,\lambda_3,\lambda_4,\lambda_5}(L,Y,Y,M,M)=(\lambda_2-\lambda_3)(\lambda_4-\lambda_5)(\lambda_2+\lambda_3+\lambda_4+\lambda_5).
\end{align*}
\end{coro}
{\it Proof.} It follows directly from the Theorem \ref{th1}, (\ref{dimc}) and (\ref{dimcc}). About the concrete expression of $\bar{\Phi},\bar{\Lambda},\bar{\Psi},\bar{\Omega},\bar{\Theta}$, we just take $\bar{\Phi}$ for example, others can be proved similarly.
\begin{align*}
\bar{\Phi}_{\lambda_1,\lambda_2}(L,L)&=(\tau(\partial\Phi))_{\lambda_1,\lambda_2}(L,L)\\
&=(-1)^2\frac{\partial}{\partial\lambda}(\partial\Phi)_{\lambda_1,\lambda_2,\lambda}(L,L,L)|_{\lambda=0}\\
&=\frac{\partial}{\partial\lambda}(\lambda_1+\lambda_2+\lambda)\Phi_{\lambda_1,\lambda_2,\lambda}(L,L,L)|_{\lambda=0}\\
&=\frac{\partial}{\partial\lambda}(\lambda_1+\lambda_2+\lambda)(\lambda_1-\lambda_2)(\lambda_1-\lambda)(\lambda_2-\lambda)|_{\lambda=0}\\
&=-\lambda_1^3+\lambda_2^3.
\end{align*}

\QED

\subsection{Computation of $H^q (\mathcal{SV},\C_a)$ if $a\neq0$}
\begin{theo} \label{th3}For any $q\in\Z_+$,
${H}^q(\mathcal{SV},\C_a)=0$ if $a\neq 0$.
\end{theo}
\noindent{\it Proof.~} Similar to the proof of Lemma 3.2 in \cite{YW}, we can define an operator $\tau_2:\t C^q(\SV,\C_a)\rightarrow \t C^{q-1}(\SV,\C_a)$ by
\begin{eqnarray}\label{7-3}
(\tau_2
\g)_{\la_1,\cdots,\la_{q-1}}(a_1,\cdots,a_{q-1})
=(-1)^{q-1}\g_{\la_1,\cdots,\la_{q-1},\la}(a_1,\cdots,a_{q-1},L)|_{\la=0}.
\end{eqnarray}
Then
\begin{align*}
(d\tau_2+\tau_2 d)
\g)\equiv -a \g \ (\mbox{mod}\
\partial\t C^q(\SV, \CC_a),
\end{align*}
which implies ${H}^q(\mathcal{SV},\C_a)=0$ for all $q\geq0$ if $a\neq 0$.
\QED

A similar argument shows
\begin{theo}\label{th4}
${H}^q(\mathcal{SV},M_{\alpha,\beta})=0$ if $\beta \neq 0$.
\end{theo}
\vskip10pt

%\noindent{ \bf{Proof of Theorem \ref{mm}.~}}
%Together with Theorem \ref{tm3.10}, \ref{tm3.11}, \ref{tm4.2} and \ref{tm4.3}, we can obtain the result immediately.

%\QED\vskip5pt
%With the similar discussion of the above two theorems, we can obtain the following result immediately.
%\begin{coro}\label{cor4.4}
%(1) ${\t H}^\bullet(\mathcal{HV},\C_a)={H}^\bullet(\mathcal{HV},\C_a)=0$ if $a\neq 0$.\\
%(2) ${\t H}^\bullet(\mathcal{HV},M_{\alpha,\beta})={H}^\bullet(\mathcal{HV},M_{\alpha,\beta})=0$ if $\beta \neq 0$.
%\end{coro}
%%%%%%%%%%%%%%%%%%5

\vspace{4mm} \noindent\bf{\footnotesize Acknowledgements.}\ \rm
{\footnotesize This work was supported by the National Natural Science Foundation of China (No. 11701345, 11871421, 12171129) and the Zhejiang Provincial Natural Science Foundation of China (No. LY20A010022) and the Fundamental Research Funds for the  Central Universities (No. 22120210554).}\\
\vskip18pt \small\footnotesize
\parskip0pt\lineskip1pt


\parskip=0pt\baselineskip=1pt
\begin{thebibliography}{9999}
\def\RE#1{\bibitem{#1}\label{#1}}


 \bibitem{BDK}
 {Barakat, A., De sole, A., Kac V.,} Poisson vertex algebras in the theory of Hamiltonian equations, Japan. J. Math. {\bf 4}, (2009) 141--252.

 \vspace{.1cm}
 \bibitem{BPZ}

{ Belavin, A., Polyakov, A., Zamolodchikov, A.,} Infinite conformal symmetry in two-dimensional quantum field theory,
 Nucl. Phys. B{\bf 241}, (1984) 333--380.
  \vspace{.1cm}


	
%%[EL]



%\bibitem{BDK} Bakalov B., D'Andrea A., Kac V., Theory of finite
%preudoalgebras, {\it Adv. Math.}, {\bf 162} (2001) 1--140.

\bibitem{BKV} Bakalov B., Kac V., Voronov A., Cohomology of
conformal algebras, Comm. Math. Phys., {\bf 200} (1999)
561--598.

\bibitem {CK} Cheng S.-J., Kac V., Conformal modules, Asian
J. Math., {\bf 1}(1) (1997) 181--193.

\bibitem {DK} D'Andrea A., Kac V., Structure theory of finite
conformal algebras, Sel. Math., New Ser., {\bf 4} (1998)
377--418.

%\bibitem {FK} Fattori D., Kac V., Classification of finite simple Lie
%conformal superalgebras, {\it J. Algebra}, {\bf 258}(1) (2002)
%23--59.

\bibitem {K1} Kac V.,  Vertex algebras for beginners, Univ. Lect. Series 10, AMS
(1996). Second edition 1998.

\bibitem {K2} Kac V., Formal distribution algebras and conformal
algebras, A talk at the Brisbane Congress in Math. Physics, July
1997, arXiv:q-alg/9709027v2.

\bibitem{KacL}
{ Kac V.,} The idea of locality, in Physical Application and Mathematical Aspects of Geometry, Groups and Algebras, edited
by H.-D. Doebner et al.et al.(World Science Publisher, Singapore,1997), pp.16-32.

   %\vspace{.1cm}

%\bibitem {K3} Kac V., The idea of locality, in: H.-D. Doebner, et al. (Eds.),
%Physical Applications and Mathematical Aspects of Geometry, Groups
%and Algebras, World Sci. Publ., Singapore, (1997) 16--32.

\bibitem{S} Su Y., Low dimensional cohomology of general conformal algebras $gc_N$,
J. Math. Phys., {\bf 45} (2004) 509--524.

\bibitem{SYX} Su Y., Yue X., Filtered Lie conformal algebras whose associated graded
algebras are isomorphic to that of general conformal algebra $gc_1$,
J. Algebra, 340 (2011) 182--198.

\bibitem{SY} Su Y., Yuan L., Schr\"{o}dinger-Virasoro Lie conformal
algebra, J. Math. Phys., {\bf54} (2013) 053503, 16pp.

\bibitem{WY} Wu H., Yuan L., Classification of finite irreducible conformal modules over some Lie conformal algebras related to the Virasoro conformal algebra, J. Math. Phys., {\bf58} (2017) 041701, 10pp.

\bibitem{YW} Yuan L., Wu H., Cohomology of the Heisenberg-Viraosro conformal algebra, J. Lie Theory, 26(2016) 1187--1197.

\bibitem{Z1} Zelmanov E., On the structure of conformal algebras, in:
Combinatorial and Computational Algebra, Hong Kong, 1999, in:
Contemp. Math., vol. 264, Amer. Math. Soc., Providence, RI, (2000)
139--153.

\bibitem{Z2} Zelmanov E., Idempotents in conformal algebras, in: Proceedings
of the Third International Algebra Conference, Tainan, 2002, Kluwer
Acad. Publ., Dordrecht, (2003) 257--266.

\end{thebibliography}
\end{document}